\documentclass[
framed_theorems,
framed_proofs,
print,
]{OCPreprint}
\usepackage{enumitem}
\setlist[enumerate,1]{label=\textup{(\roman*)}}

\let\cite\parencite

\renewcommand\d{\mathop{}\!\mathrm{d}}

\let\seq\set

\title[Newton's problem of least resistance]%
{New numerical solutions to Newton's problem of least resistance via a convex hull approach}
\author[Wachsmuth]{%
	Gerd Wachsmuth%
	\footnote{%
		Brandenburgische Technische Universität Cottbus--Senftenberg,
		Institute of Mathematics,
		03046 Cottbus,
		Germany,
		\email{gerd.wachsmuth@b-tu.de},
		\url{https://www.b-tu.de/fg-optimale-steuerung},
		ORCID: \href{https://orcid.org/0000-0002-3098-1503}{0000-0002-3098-1503}%
	}
}

\dedication{Dedicated to Alexander Plakhov on the occasion of his 65th birthday}

\begin{document}
%%fakesection: Title, Abstract, Keywords and MSC
\maketitle

\begin{abstract}
	We present a numerical method for the solution of Newton's problem of least resistance
	in the class of convex functions using a convex hull approach.
	We observe that the numerically computed solutions
	possess some symmetry.
	Further,
	their extremal points lie on several curves.
	By exploiting this conjectured structure,
	we are able to compute highly accurate solutions to Newton's problem.
\end{abstract}

\begin{keywords}
	Newton's problem of minimal resistance,
	convex hull,
	calculus of variations with convexity constraint
\end{keywords}

\begin{msc}
	\mscLink{49Q10},
	\mscLink{52A15},
	\mscLink{52A40}
\end{msc}

\section{Introduction}
\label{sec:intro}
One of the oldest problems from calculus of variations
is Newton's problem of least resistance.
In this problem,
the shape of a body which is travelling through a rare medium
has to be designed in order to minimize its resistance.
The modelling leads to the objective functional
\begin{equation*}
	J(u)
	=
	\int_\Omega \frac{1}{\abs{\nabla u}^2+1} \d x,
\end{equation*}
see \cite{ButtazzoKawohl1993}.
The functional $J$
has to be minimized over functions $u \colon \Omega \to [0,M]$,
where $\Omega \subset \R^2$ is the base of the body and $M > 0$ is its height.
In order to obtain solvability of this problem,
one further has to restrict the class of admissible functions,
see again \cite{ButtazzoKawohl1993}
or
\cite{Buttazzo2025} and the references therein.

In what follows, we consider the ``traditional'' case,
i.e., $\Omega = B_1(0) \subset \R^2$ is the (open) unit disc.
In this case,
Newton himself solved the problem in the class of concave and rotationally symmetric functions.
Almost 300 years later, it has been realized
in \cite{BrockFeroneKawohl1996}
that Newton's solution is no longer the solution in the set of concave functions
or, speaking differently, the solutions in the set of concave functions are no longer rotationally symmetric.
Since then, the problem of characterizing the optimal shapes is an open problem.
Progress has been made on the analytical side by providing necessary optimality conditions
(which have to be satisfied by the solution)
and on the numerical side by the computation of possible minimizers.

In the last decades,
many analytical properties of the minimizers have been derived.
It was shown that
$\abs{\nabla u} \not\in (0,1)$, see \cite[Theorem~2.3]{ButtazzoFeroneKawohl1995}.
Further,
if $u$ is $C^2$ on some open $\omega \subset \Omega$,
then the smallest eigenvalue of the Hessian $\nabla^2 u$ vanishes on $\omega$,
see \cite[Remark~3.4]{BrockFeroneKawohl1996}.

Further progress has been recently obtained by Alexander Plakhov.
In \cite[Theorem~2]{Plakhov2020}
he settled the open problem concerning the boundary values of solutions $u$
by showing that $u(x) \to 0$ if $x$ approaches the boundary.
He also studied the singular points of $u$,
i.e., the points $x \in \Omega$ at which $\nabla u(x)$ does not exist.
In the works \cite{Plakhov2021,Plakhov2022:1}
it was shown that the solution $u$ can be reconstructed via its singular points.
This implies that if $\nabla u(x)$ exists at all points $x$ from an open set $\omega \subset \Omega$,
then the graph of $u$ above $\omega$ does not contain extremal points of
the associated convex body $\cl\set{(x,z) \in \Omega \times \R \given 0 \le z \le u(x) }$.
Note that this generalizes the 
result from \cite[Remark~3.4]{BrockFeroneKawohl1996} mentioned above.
Finally, \cite{Plakhov2022:2} analyzes the structure of the singular points
which are so-called ridge points.

The literature concerning the numerical solution of Newton's problem is rather scarce,
the only contributions we are aware of are
\cite{Lachand-RobertOudet2005,Wachsmuth2013:1}.
Further,
there are several suggestions for (semi)-analytical descriptions of optimal bodies,
see, e.g.,
\cite{Lachand-RobertPeletier2001:1,Wachsmuth2013:1,LokutsievskiyZelikin2}.
However,
in \cite{LokutsievskiyWachsmuthZelikin2020},
it has been shown that most of these conjectured solutions
contain certain conical parts in their boundary
and, consequently, they cannot be optimal.
In particular,
the cases in which $M$ is bigger than roughly $1.5$
and in which $M \in [0.95, 1.08]$
remain open,
see \cite[Section~6]{LokutsievskiyWachsmuthZelikin2020}.

We also mention \cite{MaininiMonteverdeOudetPercivale2019}.
Therein, a variation of the problem
(concavity of $u$ is replaced by semiconcavity)
is studied and numerical solutions are presented.

The goal of the present paper
is to provide numerical solutions of higher accuracy.
We approximate the graph of the function $u$ as
the convex hull of finitely many points.
The position of these points are optimized
in order to minimize $J(u)$,
see \cref{sec:cvx} for details.
The numerical results
presented in \cref{sec:free}
suggest that the optimal bodies
can be represented as the convex hull of finitely many extremal arcs,
see in particular \cref{fig:free_solutions}.
Thus, it seems worthwhile to exploit this conjectured structure of the solutions.
This leads to the method described in \cref{sec:restricted}.
With this improvement, it is possible to obtain highly accurate approximations,
see \cref{fig:restricted_solutions}.
We can also confirm the non-optimality result from \cite{LokutsievskiyWachsmuthZelikin2020},
see \cref{tab:conjectured_solutions_detail}.

\section{Convex hull approach}
\label{sec:cvx}

We suggest to approximate (or discretize) the graph of the concave function $u \colon \Omega \to [0,M]$
as the (upper) boundary of the convex hull of finitely many points.
To be precise,
given finitely many points $\seq{x_i}_{i = 1,\ldots,n} \subset \Omega \times (0,M]$,
we consider
the convex body
\begin{equation*}
	P(\seq{x_i}_i)
	=
	\conv\parens*{
		(\cl\Omega \times \set{0}) \cup \set{ x_i \given i = 1,\ldots, n}
	}
	.
\end{equation*}
Then, the resistance of this body can be computed via
\begin{equation}
	\label{eq:formula_objective}
	f(\seq{x_i}_i)
	:=
	\int_{\partial^+ P(\seq{x_i}_i)} \nu_3 \d\HH^2,
\end{equation}
where $\partial^+ P(\seq{x_i}_i)$ is the upper boundary of the body, i.e.,
\begin{equation*}
	\partial^+ P(\seq{x_i}_i)
	=
	\partial P(\seq{x_i}_i)
	\setminus (\Omega \times \set{0}),
\end{equation*}
$\nu_3$ is the third component of the outer unit normal vector
and $\HH^2$ is the two-dimensional Hausdorff measure,
see \cite[Remark~3.2]{Buttazzo2025}.

In case that the points $\seq{x_i}_i$ are in general position,
it is easy to see that the upper boundary of
$P(\seq{x_i}_i)$
consists of finitely many parts of the following types
\begin{itemize}
	\item
		triangles spanned by three distinct points $x_i, x_j, x_l$,
	\item
		triangles spanned by two points $x_i, x_j$ and a point $y_l \in \partial\Omega \times \set{0}$,
	\item
		conical pieces which are the convex hull of a point $x_i$ and a circular segment $C_l \subset \partial\Omega \times \set{0}$.
\end{itemize}
Now, it is almost straightforward to evaluate the integral in \eqref{eq:formula_objective}.
For both types of triangles, one just has to evaluate the (constant) normal vector on this triangle.
For the conical pieces, one arrives at an integral of type
\begin{equation*}
	\int_a^b \frac{(1 - \mu\cos(t))^3}{M^2 + (1 - \mu\cos(t))^2} \d t
\end{equation*}
with a parameter $\mu > 0$,
which can be computed via the approach given in \cite[Appendix~A]{Lachand-RobertPeletier2001:1}.

Further,
small changes of the coordinates $\seq{x_i}_i$ do not change the structure of
the boundary of the body $P(\seq{x_i}_i)$.
Consequently, the function $f$ is smooth
at all points $\seq{x_i}_i$ which are in general position
and it is straightforward to compute first and second derivatives in these points.
The evaluation of the function $f$ readily extends to points
which are not in general position.
In this case, one could have polygons instead of the triangles mentioned above.
It is, however, not clear if the function $f$ is still differentiable at those points.

Although the solutions to Newton's problem are not radially symmetric,
all the numerical solutions reported in the literature possess symmetry
given by the dihedral group $D_k$, $k \ge 2$,
which is the symmetry group of the regular $k$-sided polygon.
Thus, it seems to be a good idea to exploit this symmetry in the computations.
Note that this is easily possible with the described convex hull approach.
To this end, let
$\Pi_k$ be the mapping which maps every point to its orbit in the group $D_k$,
i.e.,
(unless the point lies on the axes of symmetry)
it is mapped to its $2 k$ images under the group elements of $D_k$.
Then, instead of using $f$ from \eqref{eq:formula_objective} directly,
one would consider
\begin{equation*}
	\label{eq:formula_objective_symmetry}
	f_k( \seq{x_i}_i )
	:=
	f(\Pi_k \seq{x_i}_i )
	.
\end{equation*}
Again, one can check that $f_k$ is differentiable outside of a negligible set.

\section{Free optimization}
\label{sec:free}
In this section,
we consider an optimization procedure without assuming any structure of the body
(besides some possible symmetry).
For a given number of points $n \in \N$,
we consider the problems
\begin{equation}
	\label{eq:free_nonsymmetric}
	\text{Minimize} \qquad
	f(\seq{x_i}_i)
	\qquad\text{with respect to}
	\qquad
	x_i \in \Omega \times (0,M] \; \forall i=1,\ldots,n
\end{equation}
for the case without symmetry,
or, given additionally a symmetry parameter $k \ge 2$,
\begin{equation}
	\label{eq:free_symmetric}
	\text{Minimize} \qquad
	f_k(\seq{x_i}_i)
	\qquad\text{with respect to}
	\qquad
	x_i \in S_k \times (0,M] \; \forall i=1,\ldots,n,
\end{equation}
where
$S_k := \set{(r \cos(\varphi), r \sin(\varphi)) \given r \in [0,1), \varphi \in [0,\pi/k]}$ is a sector with opening angle $\pi / k$.

We briefly mention that we combined a projected gradient method
with a regularized Newton method on the inactive components
for the numerical solution of these problems.
Note that these problems are highly nonconvex
and we can only expect that the numerical methods compute
local minimizers or stationary points.
After each run of the optimization algorithm,
we perturbed the solution by adding additional points
which lie slightly above the upper boundary of the convex body.
In this way, we are able to refine the solution
and the optimization method can be started again.

\begin{figure}[htp]
	\centering
	\hfill%
	\includegraphics[width=.45\textwidth]{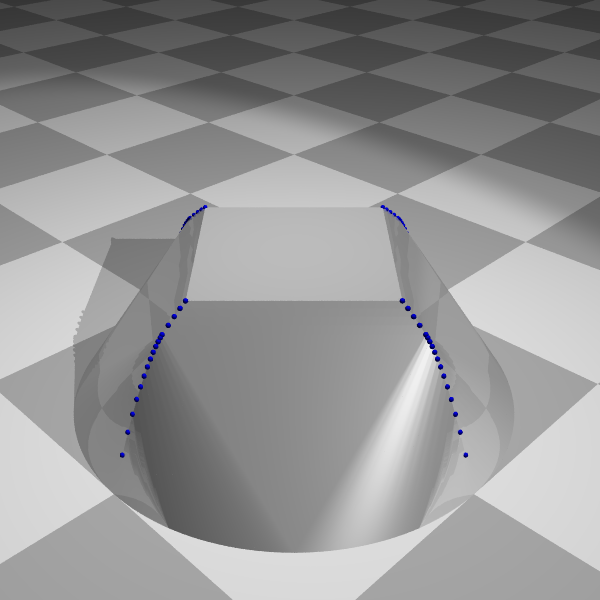}%
	\hfill%
	\includegraphics[width=.45\textwidth]{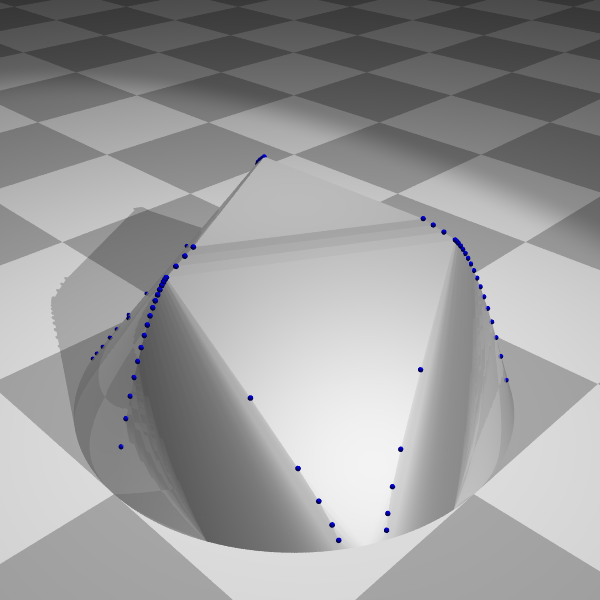}%
	\hfill%
	\\[.5cm]
	\hfill%
	\includegraphics[width=.45\textwidth]{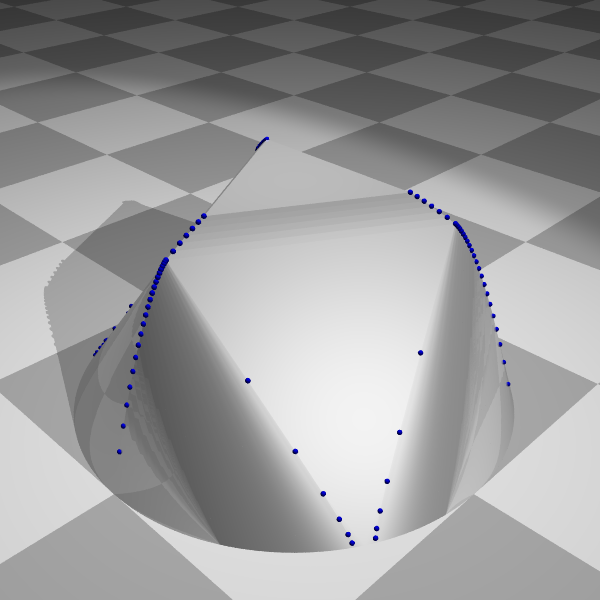}%
	\hfill%
	\includegraphics[width=.45\textwidth]{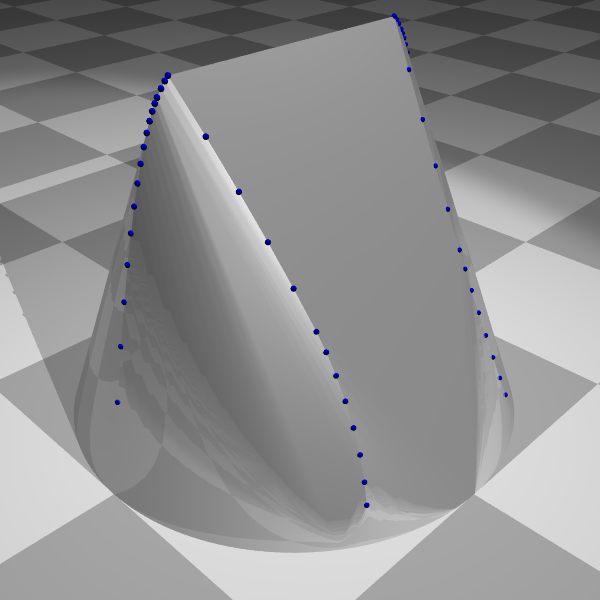}%
	\hfill%
	\caption{%
		Numerically obtained solutions of \eqref{eq:free_symmetric}
		for the parameters
		$M = 0.7$, $k = 4$ (top left),
		$M = 0.9$, $k = 3$ (top right),
		$M = 1.0$, $k = 3$ (bottom left),
		$M = 1.5$, $k = 2$ (bottom right).
	}
	\label{fig:free_solutions}
\end{figure}

Let us report some of the numerical results.
The solution of \eqref{eq:free_nonsymmetric},
in which the symmetry was not prescribed,
always resulted in bodies which possessed almost the symmetry given by $D_k$
for some value of $k \ge 2$ (depending on the height parameter $M > 0$).
Consequently,
it is beneficial to consider the symmetric version \eqref{eq:free_symmetric} instead
for two reasons.
First,
the symmetry can be used to reduce the dimension of the problem by a factor of
roughly $(2m)^{-1}$.
Second,
the objective in \eqref{eq:free_nonsymmetric} is rotationally invariant
which means that its Hessian cannot be positive definite
and this is cumbersome for optimization methods.
Note that this rotational invariance is removed in \eqref{eq:free_symmetric}.
Therefore,
we stick to the solution of \eqref{eq:free_symmetric} from now on.
We solved \eqref{eq:free_symmetric}
for four combinations of $M$ and $k$,
see \cref{fig:free_solutions} and \cref{tab:free_solutions}.
\begin{table}
	\centering
	\begin{tabular}{ccccc}
		\toprule
		$M$ & $k$ & $n$ & objective value & runtime [s] \\
		\midrule
		$0.7$ & $4$ & $17$ & $1.5454863290500$ & $ 8$ \\
		$0.9$ & $3$ & $25$ & $1.2620157197502$ & $15$ \\
		$1.0$ & $3$ & $34$ & $1.1377471863594$ & $29$ \\
		$1.5$ & $2$ & $29$ & $0.6999056973373$ & $12$ \\
	\bottomrule
	\end{tabular}
	\caption{Numerical results for \eqref{eq:free_symmetric}.}
	\label{tab:free_solutions}
\end{table}%
In \cref{fig:free_solutions},
we showed the location of the points $\Pi_k \seq{x_i}_i$
by small, blue dots.
In the case $M = 0.7$, $k = 4$,
we see that these points only lie on the symmetry axes.
Thus, this computation confirms
that the structural conjecture from \cite[Section~3]{Wachsmuth2013:1}
is satisfied for this choice of $M$.
We recall that the non-optimality from \cite{LokutsievskiyWachsmuthZelikin2020}
does not apply in this case.
In the other cases,
some of the points $\seq{x_i}_i$ do not lie on the axes of symmetry.
Instead, it seems that they are forming a second curve of extremal points
which, due to the application of $\Pi_k$, are copied to $2 k$ curves of extremal points.
Note that this is in line with the observation of \cite{LokutsievskiyWachsmuthZelikin2020},
see in particular Figure~3 therein.
In the discussed cases,
if we build the solution just as the convex hull of curves lying in the symmetry planes,
it contains conical parts
which are provably not optimal.
In the solutions shown in \cref{fig:free_solutions},
these conical parts are broken up by the additional points
which do not lie in the planes of symmetry.

In \cref{tab:free_solutions} we show some numbers for the solution of \eqref{eq:free_symmetric}.
The reported run time was obtained using
Matlab R2024b
on a
Intel Core i7-8700.
We see that this approach is able
to compute some nice approximations in rather short time.
By comparing the objective values in \cref{tab:free_solutions}
with the more accurate values presented in \cref{tab:restricted_solutions}
below,
we see that 3 to 4 leading digits of the objective value
are already accurate.
Moreover,
the objective values improve the results given in \cite[Table~2]{Wachsmuth2013:1}
although the number of unknowns and the computational time
are much smaller
in our new approach.

\section{Restricted optimization}
\label{sec:restricted}
We have seen in \cref{fig:free_solutions},
that the numerically obtained solutions
seem to possess a very distinctive structure.
The extremal points (small blue dots)
lie on two curves
(before the application of the symmetric copying operator $\Pi_k$):
One curve lies directly on the plane of symmetry and ends at the top face of the body
and a second one lies between the planes of symmetry.
The second curve is not present in the case $M = 0.7$.
It goes to the top face in the case $M = 1.5$
and in the other two cases, it ends on the first curve of extremal points
below the top face of the body.

It seems worthwhile to exploit this particular structure of the solutions.
We are doing this in the following way.
Given $n_1, n_2 \in \N$,
we consider the optimization variables
$z \in (0,1)$,
$\seq{Y_i}_{i = 1,\ldots,n_1} \subset (0,1) \times (0,M)$
and
$\seq{X_i}_{i = 1,\ldots,n_2} \subset \Omega \times (0,M)$.
The meaning of these variables is as follows.
The scalar $z$ describes where the first extremal curve on the symmetry plane $\set{x_1 = 0}$
hits the upper face of the convex body,
i.e., $z$ corresponds to the point $(0, z, M)$.
The vectors $\seq{Y_i}_{i = 1,\ldots,n_1}$ describe the extremal curve on the symmetry plane $\set{x_1 = 0}$,
i.e., the vector $Y_i \in \R^2$ corresponds to the point $(0, Y_{i,1}, Y_{i,2})$.
Finally, the vectors $\seq{X_i}_{i = 1,\ldots,n_2}$ describe the second extremal curve
which does not lie in the symmetry plane $\set{x_1 = 0}$,
i.e., the vector $X_i \in \R^3$ corresponds directly to the point $X_i$.

This being said, the overall objective is now given by
\begin{equation*}
	f_{k, \textup{restricted}}(z, \seq{Y_i}, \seq{X_i})
	=
	f_k\parens*{
		\seq*{
			(0,z,M),
			(0, Y_{i,1}, Y_{i,2}),
			X_j
			\given
			\begin{aligned}
				& i = 1,\ldots, n_1, \\
				& j = 1,\ldots, n_2
			\end{aligned}
		}
	}
	.
\end{equation*}
Consequently,
we arrive at the optimization problem
\begin{equation}
	\label{eq:restricted_symmetry}
	\begin{aligned}
		\text{Minimize} \quad &f_{k, \textup{restricted}}(z, \seq{Y_i}, \seq{X_i}) \\
		\text{with respect to} \quad &
		z \in (0,1), \;
		\seq{Y_i}_{i = 1,\ldots,n_1} \subset (0,1) \times (0,M),
		\\&
		\seq{X_i}_{i = 1,\ldots,n_2} \subset \Omega \times (0,M).
	\end{aligned}
\end{equation}
Note that the feasible set is an open set.
Although we cannot prove it,
it is very reasonable that the objective still possesses minimizers over the open, feasible set:
\begin{itemize}
	\item
		$z$ cannot go to $0$, since the optimal bodies possess a flat top face,
	\item
		$z$ and $Y_{i,1}$ cannot go to $1$, since optimal functions vanish on the boundary, \cite[Theorem~2]{Plakhov2020};
		for the same reason, $(X_{i,1}, X_{i,2})$ has to stay inside $\Omega$,
	\item
		$Y_{i,2}$ and $X_{i,3}$ stay strictly above $0$, since otherwise the corresponding point
		would not be visible in the convex hull,
	\item
		$Y_{i,2}$ and $X_{i,3}$ stay strictly below $M$, since the top face is determined by the value of $z$.
\end{itemize}
Consequently,
we can use an optimization method for unconstrained optimization problems
to minimize $f_{k, \textup{restricted}}$
over its open, feasible set.
We only have to bound the steps taken in the optimization method
such that the iterates do not leave the open, feasible set
due to the taken steps being too large.
In our implementation,
we used a regularized Newton method.

\begin{figure}[htp]
	\centering
	\hfill%
	\includegraphics[width=.45\textwidth]{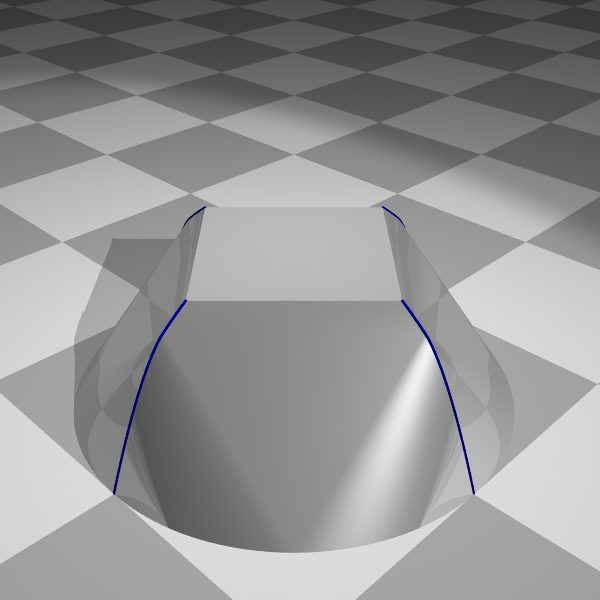}%
	\hfill%
	\includegraphics[width=.45\textwidth]{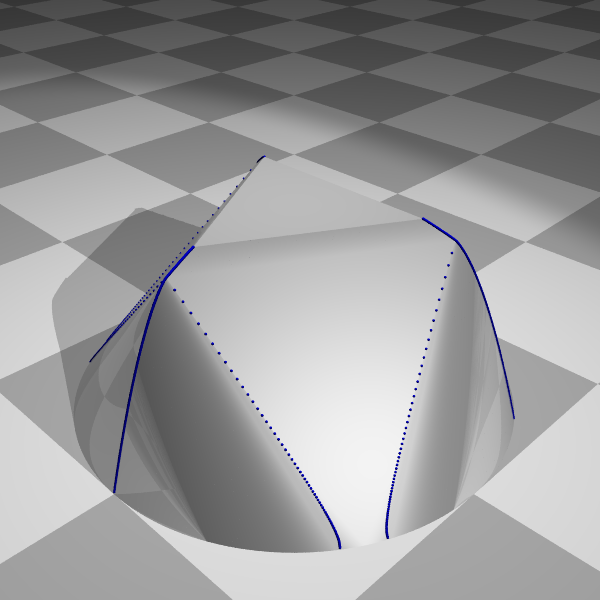}%
	\hfill%
	\\[.5cm]
	\hfill%
	\includegraphics[width=.45\textwidth]{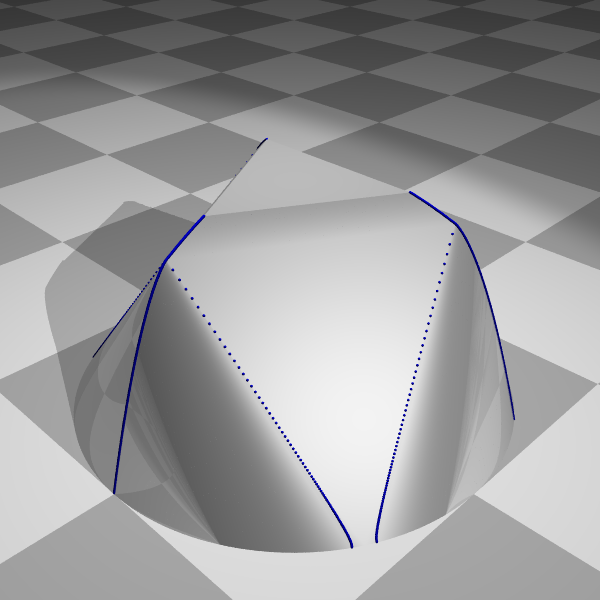}%
	\hfill%
	\includegraphics[width=.45\textwidth]{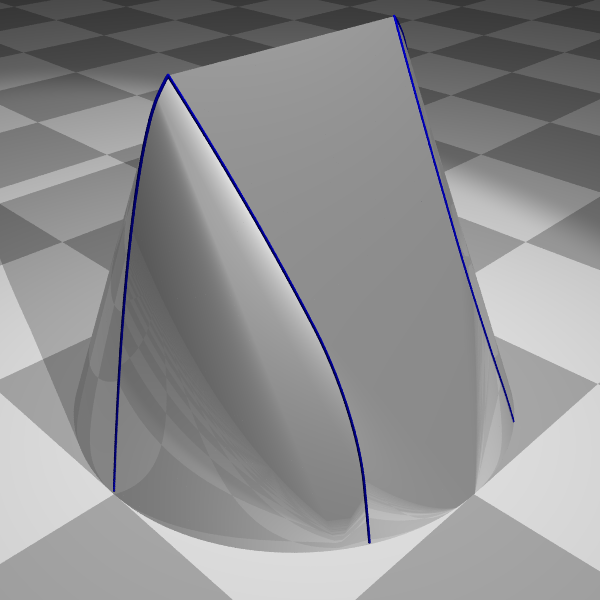}%
	\hfill%
	\caption{%
		Numerically obtained solutions of \eqref{eq:restricted_symmetry}
		for the parameters
		$M = 0.7$, $k = 4$ (top left),
		$M = 0.9$, $k = 3$ (top right),
		$M = 1.0$, $k = 3$ (bottom left),
		$M = 1.5$, $k = 2$ (bottom right).
	}
	\label{fig:restricted_solutions}
\end{figure}

After a potential minimizer of 
$f_{k, \textup{restricted}}$
with given numbers $n_1, n_2$
has been found,
it is possible
construct an initial point for higher values of $n_1$ or $n_2$
by a refinement strategy.
For example, $\seq{Y_i}$ could be replaced by
(assuming that $\seq{Y_{i,1}}$ is ordered increasingly)
\begin{equation*}
	\seq*{
		\begin{aligned}
			&
			( (z + Y_{i,1})/2, (M + Y_{i,2}) / 2 + \varepsilon),
			(Y_{1,1}, Y_{1,2}),
			( (Y_{1,1} + Y_{2,1})/2, (Y_{1,2} + Y_{2,2}) / 2 + \varepsilon),
			\\ &
			(Y_{2,1}, Y_{2,2}),
			( (Y_{2,1} + Y_{3,1})/2, (Y_{2,2} + Y_{3,2}) / 2 + \varepsilon),
			(Y_{3,1}, Y_{3,2}),
			\ldots,
			(Y_{n_1,1}, Y_{n_1,2}),
			\\ &
			((Y_{n_1,1} + 1)/2, (Y_{n_1,2} + 0)/2 + \varepsilon)
		\end{aligned}
	}
	,
\end{equation*}
i.e., we have inserted interpolations of neighboring points into $\seq{Y_i}$.
Note that $\varepsilon > 0$ is necessary such that the newly inserted points
are visible in the convex hull.
A similar idea can be used to refine the variable $\seq{X_i}$.
With this approach,
it is possible
to compute highly accurate solutions
to Newton's problem of least resistance.
Our results are shown in \cref{fig:restricted_solutions}
and \cref{tab:restricted_solutions}.
\begin{table}[htp]
	\centering
	\begin{tabular}{cccccc}
		\toprule
		$M$ & $k$ & $n_1$ & $n_2$ & objective value & runtime [s] \\
		\midrule
		$0.7$ & $4$ & $8556$ & $   0$ & $1.545460647$ & $458$ \\
		$0.9$ & $3$ & $6062$ & $ 256$ & $1.261987252$ & $247$ \\
		$1.0$ & $3$ & $5643$ & $  96$ & $1.137729313$ & $354$ \\
		$1.5$ & $2$ & $2972$ & $1986$ & $0.699881606$ & $ 67$ \\
		\bottomrule
	\end{tabular}
	\caption{Numerical results for \eqref{eq:restricted_symmetry}.}
	\label{tab:restricted_solutions}
\end{table}%
It can be seen,
that our approach is capable of
resolving the extremal arcs.
The runtime is still acceptable
and compared with \cref{tab:free_solutions},
we can further improve the computed objective values.

Finally, we make some computations complementing
\cite[Table~2]{LokutsievskiyWachsmuthZelikin2020}.
Therein, it has been shown that the conjectured structure from
\cite[Section~3]{Wachsmuth2013:1}
for heights $M \in [0.95, 1.08]$
cannot be optimal since the corresponding minimizers contain nonoptimal conical parts.
For $21$ values between $M = 0.9$ and $M = 1.1$,
we applied our algorithm to improve
the objective values given in
\cite[Table~2]{LokutsievskiyWachsmuthZelikin2020}.
The results are shown in \cref{tab:conjectured_solutions_detail}.
\begin{table}[htp]
	\centering
	\begin{tabular}{*{4}l}
		\toprule
		& \multicolumn{2}{c}{Structured solution} & New results \\
		\cmidrule(r){2-3} 
		\cmidrule(r){4-4} 
		$M$  & $k = 3$ & $k = 4$ & $k = 3$\\
		\midrule
		$0.90$ &             $1.2619895$ CN  & \underline{$1.2599052$} &            $1.2619873$  \\
		$0.91$ &             $1.2488725$ CN  & \underline{$1.2472919$} &            $1.2488707$  \\
		$0.92$ &             $1.2359138$ CN  & \underline{$1.2348243$} &            $1.2359124$  \\
		$0.93$ &             $1.2231117$ CN  & \underline{$1.2225019$} &            $1.2231106$  \\
		$0.94$ &             $1.2104642$ CN  & \underline{$1.2103218$} &            $1.2104634$  \\
		$0.95$ &             $1.1979697$ CN  &            $1.1982830$  & \underline{$1.1979691$} \\
		$0.96$ &             $1.1856264$ CN  &            $1.1863837$  & \underline{$1.1856259$} \\
		$0.97$ &             $1.1734324$ CN  &            $1.1746215$  & \underline{$1.1734320$} \\
		$0.98$ &             $1.1613861$ CN  &            $1.1629969$  & \underline{$1.1613858$} \\
		$0.99$ &             $1.1494856$ CN  &            $1.1515072$  & \underline{$1.1494855$} \\
		$1.00$ &             $1.1377294$ CN  &            $1.1401510$  & \underline{$1.1377293$} \\
		$1.01$ &             $1.1261157$ CN  &            $1.1289274$  & \underline{$1.1261156$} \\
		$1.02$ &             $1.1146427$ CN  &            $1.1178333$  & \underline{$1.1146427$} \\
		$1.03$ &             $1.1033089$ CN  &            $1.1068681$  & \underline{$1.1033089$} \\
		$1.04$ &             $1.0921125$ CN  &            $1.0960303$  & \underline{$1.0921125$} \\
		$1.05$ &             $1.0810520$ CN  &            $1.0853184$  & \underline{$1.0810520$} \\
		$1.06$ &             $1.0701256$ CN  &            $1.0747311$  & \underline{$1.0701256$} \\
		$1.07$ &             $1.0593317$ CN  &            $1.0642667$  & \underline{$1.0593317$} \\
		$1.08$ &             $1.0486688$ CN  &            $1.0539240$  & \underline{$1.0486688$} \\
		$1.09$ &  \underline{$1.0381352$ C } &            $1.0437014$  & \underline{$1.0381352$} \\
		$1.10$ &  \underline{$1.0277294$   } &            $1.0335975$  & \underline{$1.0277294$} \\
	\end{tabular}
	\caption{%
		We display computed function values for heights $M \in [0.9, 1.1]$.
		The first two columns
		(``Structured solution'')
		are from \cite[Table~2]{LokutsievskiyWachsmuthZelikin2020}
		in which we evaluated the
		conjectured optimal values using the conjecture from \cite[Section~3]{Wachsmuth2013:1}.
		A ``C'' indicates that this solution contains a conical part.
		For the solutions marked by ``N'' the conical part is nonoptimal according to \cite[Corollary~7]{LokutsievskiyWachsmuthZelikin2020}.
		The final column shows the computational results using \eqref{eq:restricted_symmetry}.
		The best known solution for each $M$ are underlined.
	}
	\label{tab:conjectured_solutions_detail}
\end{table}
We can observe
that we are able to improve all solutions
in which a nonoptimal conical part is present.
However, the removal of the conical part
improves the objective value only very slightly.
In fact, for $M \in [1.02, 1.08]$,
the nonoptimality of the conical part is not visible
in the displayed digits.
For smaller values $M \in [0.9, 1.01]$,
we see a very slight improvement.
We finally mention that the
computational time for each row of \cref{tab:conjectured_solutions_detail}
was between one and four minutes.

\section{Summary}
We used a structured discretization approach to address Newton's problem of least resistance
in the class of convex bodies.
The obtained solutions improve other solutions which are reported in the literature.
For height values $M \in [0.95, 1.08] \cup [1.5, \infty)$
we conjecture that the optimal body
can be described via two extremal arcs
one if which is located in the plane of symmetry, see \cref{fig:restricted_solutions}.
Note that this is a generalization of the structure described in
\cite[Section~3]{Wachsmuth2013:1}.
We think that it should be possible to exploit this structure
and to obtain ordinary differential equations satisfied by the extremal arcs,
cf.\ \cite[Section~3]{Wachsmuth2013:1}.
This is subject to future research.

%%fakesection: Bib
\printbibliography

\end{document}